\documentclass[12pt]{amsart}
\usepackage{amsmath,amsthm,amsfonts,amssymb}
\usepackage{graphicx}
\usepackage{mathrsfs}
\usepackage{enumitem}
\usepackage{etoolbox}
\usepackage{xcolor}
\apptocmd{\sloppy}{\hbadness 10000\relax}{}{}
\apptocmd{\sloppy}{\vbadness 10000\relax}{}{}
\usepackage[pdfpagelabels]{hyperref}
\usepackage[letterpaper,margin=1.1in]{geometry}

\numberwithin{equation}{section}
\theoremstyle{plain}

\newtheorem{theorem}{Theorem}[section]

\newtheorem{corollary}[theorem]{Corollary}

\newtheorem{lemma}[theorem]{Lemma}

\theoremstyle{definition}
\newtheorem{remark}[theorem]{Remark}
\newtheorem{definition}[theorem]{Definition}

\newcommand{\norm}[1]{\left\lVert#1\right\rVert}

\def\XXint#1#2#3{{\setbox0=\hbox{$#1{#2#3}{\int}$ }
\vcenter{\hbox{$#2#3$ }}\kern-.6\wd0}}

\newcommand{\dist}{\mathop\mathrm{dist}\nolimits}

\title{Non-existence of cusps for degenerate Alt-Caffarelli functionals}
\author{Sean McCurdy and Lisa Naples}


\address{Department of Mathematics\\ National Taiwan Normal University\\ Taipei, Taiwan}
\email{smccurdy@ntnu.edu.tw}

\address{Department of Mathematics, Statistics, and Computer Science\\ Macalester College\\1600 Grand Avenue\\Saint Paul, MN 55105}
\email{lnaples@macalester.edu}

\begin{document}

\maketitle

\begin{abstract} We eliminate the existence of cusps in a class of \textit{degenerate} free-boundary problems for the Alt-Caffarelli functional
$J_{Q}(v, \Omega):= \int_{\Omega}|\nabla v|^2 + Q^2(x)\chi_{\{v>0\}}dx,$ so-called because $Q(x) = \dist(x, \Gamma)^{\gamma}$ for $\Gamma$ an affine $k$-plane and $0< \gamma$. This problem is inspired by a generalization of the variational formulation of the Stokes Wave by \cite{AramaLeoni12}.  The elimination of cusps implies that the results of \cite{Mccurdy20} in fact describe the entire free-boundary as it intersects $\Gamma$.
\end{abstract}

\section{Introduction}

In this note, we continue to investigate the geometry of free-boundaries arising from local minimizers of a broad class of \emph{degenerate} Alt-Caffarelli functionals
\begin{align}\label{e: AC functional}
    J_{Q}(v, \Omega):= \int_{\Omega}|\nabla v|^2 + Q^2(x)\chi_{\{v>0\}}dx
\end{align}
where $\Omega\subset \mathbb{R}^n$ is an open set with Lipschitz boundary, $Q:\Omega \rightarrow \mathbb{R}_+ $ is continuous, and $n \ge 2$. A function $u$ is a \textit{minimizer} of (\ref{e: AC functional})
in the class $$K_{u_0, \Omega}: = \{u \in W^{1, 2}(\Omega) : u-u_0 \in W^{1,2}_0(\Omega) \}$$ for a $u_0 \in W^{1, 2}(\Omega)$ satisfying $u_0 \ge 0$ if for every function $v \in K_{u_0, \Omega},$ $J_{Q}(u, \Omega) \le J_{Q}(v, \Omega)$.  A function $u$ is called an \textit{$\epsilon_0$-local minimizer} of $J_Q(\cdot, \Omega)$ if there exists an $0<\epsilon_0$ such that $J_{Q}(u, B_r(x)) \le J_{Q}(v, B_r(x))$ for every $v \in K_{u, \Omega}$ satisfying
\begin{align}
\norm{\nabla (u-v)}^2_{L^2} + \norm{\chi_{\{u>0\}} - \chi_{\{v>0\}}}_{L^1(\Omega)} < \epsilon_0.
\end{align}
For any function $u$ which is an $\epsilon_0$-local minimizer of (\ref{e: AC functional}) the set $\partial \{u>0\} \cap \Omega$ is called the \emph{free-boundary}. 

The geometry of free-boundaries have been the subject of keen interest.  In particular, the \emph{non-degenerate} case where we assume,
$$
0< Q_{min} \le Q(x) \le Q_{max}< \infty,
$$
for all $x \in \Omega$, has been a source of much work (see for example \cite{AltCaffarelli81}, \cite{Weiss99}, \cite{CaffarelliJerisonKenig04},  \cite{JerisonSavin15}, \cite{DeSilvaJerison09}, \cite{EdelenEngelstein19}). Without attempting to summarize the known results in the non-degenerate case, we focus our attention on the roles that $Q_{min}, Q_{max}$ play in the weak geometry of the free-boundary. Classically, the assumption that $0<Q_{min}$ and $Q_{\max}<\infty$ imply that $\epsilon_0$-local minimizers $u$ satisfy weak geometric regularity (interior and exterior ball condition) at sufficiently small scales depending upon $\epsilon_0$ \cite{AltCaffarelli81} (see  \cite{DavidToro15} for similar results on almost-minimizers). This directly implies that in the \textit{non-degenerate} case all $x \in \partial \{u>0\} \cap \Omega$ satisfy
\begin{align*}
\frac{\mathscr{H}^n(\{u>0\} \cap B_r(x))}{\omega_n r^n} \in [c(n, Q_{min}), C(n, Q_{max})] \in (0, 1)
\end{align*}
for all $r \in (0, \min\{r(\epsilon_0), \frac{1}{2}\dist(x, \partial \Omega)\})$. In particular, the density of the positivity set $\{u>0\}$ satisfies $\Theta^n_{\{u>0\}}(x)\in [c(n, Q_{min}), C(n, Q_{max})].$ 

However, in the \emph{degenerate} case, where $0=Q_{min},$ the techniques used to prove $$\frac{\mathscr{H}^n(\{u>0\} \cap B_r(x))}{\omega_n r^n} \ge c > 0$$ fail in balls $B_r(x)$ where $Q_{min}=0$. Therefore, one of the most fundamental questions in the theory of the geometry of the free-boundary in the \emph{degenerate} case is whether the cusp set
\begin{align*}
    \Sigma := \{x \in \partial \{u>0\}: \Theta^n_{\{u>0\}}(x) = 0 \}
\end{align*}
is non-empty or not.  Investigation of the degenerate case was inaugurated by \cite{AramaLeoni12} under the assumption that $n=2$, $Q(x, y) = \sqrt{(h-y)_+}$ and $\Omega = [0, 1] \times [0, \infty)$ and $0<h<\infty$, with subsequent generalization to higher dimensions and other exponents by \cite{GravinaLeoni18, GravinaLeoni19}.  Because the investigations of \cite{AramaLeoni12,GravinaLeoni18,GravinaLeoni19} were inspired by variational models of the Stokes Wave, they made strong assumptions of symmetry which allow them to trivially eliminate cusps. Inspired by this previous work the first author studied the cases $n \ge 2$ and $Q(x) = \dist(x, \Gamma)^{\gamma}$ for $\Gamma$ a $C^{1, \alpha}$-submanifold of dimension $0 \le k \le n-1$ and $0< \gamma$ without assumptions of symmetry \cite{Mccurdy20}. Among other things, this work proved that $\mathcal{H}^{n-1}(\Sigma) = 0$, but little else was able to be said about $\Sigma$. In a recent follow-up paper \cite{McCurdy21}, the first author proved that for $n=2$, $Q(x,y) = |y|^{\gamma}$ for $0< \gamma$ and assuming that $\{u>0 \} \cap \{(x, 0): x \in \mathbb{R}\} = \emptyset$, then $\Sigma = \emptyset$. While this case is sufficient to cover the two-dimensional variational formulation of the Stokes Wave in \cite{AramaLeoni12} (with symmetry assumptions removed), it left open the question of cusps in higher dimensions ($n > 2$), higher co-dimensions ($0 \le k \le n-1$), and the case that $\{u>0 \} \cap \{(x, 0)\in \mathbb{R}^k \times \mathbb{R}^{n-k}: x \in \mathbb{R}^k\} \not = \emptyset$. The focus of this note is to eliminate cusps in these cases.  

We continue to work within the framework of \cite{McCurdy21}, assuming that $\Gamma = \{(x, 0) \in \mathbb{R}^k \times \mathbb{R}^{n-k}\}$ and $Q(x, y) = |y|^\gamma$ for some $0<\gamma$.  We write $\mathbb{R}^n = \mathbb{R}^k \times \mathbb{R}^{n-k}$ and denote $B^k_r(x) \subset \mathbb{R}^k$ the ball of radius $r$ around the point $x \in \mathbb{R}^k$. With this notation we prove the following theorem.

\begin{theorem}(Main Theorem)\label{t:main theorem}
Let $0< \epsilon_0$, and let $n, k$ be integers such that $n \ge 2$ and $0 \le k \le n-1$. Let $\Gamma = \{(x, 0) \in \mathbb{R}^k\times \mathbb{R}^{n-k}\}$, and $0< \gamma$. Let $Q(x, y)= |y|^{\gamma}$. If $u$ is an $\epsilon_0$-local minimizer of $J_Q(\cdot, B_2^n(0, 0))$, then $\Sigma = \emptyset.$
\end{theorem}

\begin{corollary}(Applying \cite{Mccurdy20})
Let $0< \epsilon_0$, and let $n, k$ be integers such that $n \ge 2$ and $0 \le k \le n-1$.  Let $(x, y) \in \mathbb{R}^{k}\times \mathbb{R}^{n-k}$ and $Q(x, y) = |y|^k.$ Let $u$ be an $\epsilon_0$-local minimizer of (\ref{e: AC functional}). Then the results of \cite{Mccurdy20} describe $\partial \{u>0\} \cap \Gamma.$ That is, the set $\partial \{u>0\} \cap \Gamma = \mathcal{S}$ is countably $\min\{k,n-2\}$-rectifiable and satisfies finite upper Minkowski content bounds
\begin{align*}
Vol(B_r(\mathcal{S}\cap B_1(0))) \le C(n, \epsilon, \gamma)r^{n-\min\{n-2, k\}}.
\end{align*}
\end{corollary}
See \cite{Mccurdy20} for further details and estimates on the quantitative strata.

In fact, Theorem \ref{t:main theorem} follows from the more general proposition concerning the geometry of the positivity set $\{u>0\}$ for $\epsilon_0$-local minimizers of (\ref{e: AC functional}) under the assumptions of Theorem \ref{t:main theorem}.

\begin{definition}(Property P)
Let $0< \epsilon_0$, and let $n, k$ be integers such that $n \ge 2$ and $0 \le k \le n-1$. Let $\Gamma = \{(x, 0) \in \mathbb{R}^k\times \mathbb{R}^{n-k}\}$, and $0< \gamma$. Let $Q(x, y)= |y|^{\gamma}$. Let $u$ be an $\epsilon_0$-local minimizer of $J_Q(\cdot, B_2^{n}(0,0))$.  For an integer $N \in \mathbb{N}$, a ball $B_{r}^{n}(x, 0) \subset B_2^n(0,0)$, is said to satisfy \emph{Property P with constant N} if
\begin{align*}
    B_{r}^{k}(x) \times B^{n-k}_{\frac{2r}{N}}(0) \cap \{u>0\}
\end{align*}
has a component $\mathcal{O}$ which satisfies the following conditions.
\begin{enumerate}
\item (Intersection) $\mathcal{O} \cap \{x\}\times [-2r/N, 2r/N] \not = \emptyset$.
\item (Upper Height Bound) $$\sup_{0 \le \rho \le r}\sup\{|y| :(x', y) \in \mathcal{O}, |x' - x| = \rho\} \le \frac{r}{N}.$$ 
    \item (Lower Height Bound) $$\inf_{0\le \rho \le r} \sup \{|y|: (x', y) \in \mathcal{O}, |x'-x| = \rho \} \ge \frac{r}{N4\cdot 2^{1+\gamma}}.$$
\end{enumerate}
\end{definition}

\begin{theorem}(Technical Theorem)\label{t:main theorem 2}
Let $0< \epsilon_0$, and let $n, k$ be integers such that $n \ge 2$ and $0 \le k \le n-1$. Let $\Gamma = \{(x, 0) \in \mathbb{R}^k\times \mathbb{R}^{n-k}\}$, and $0< \gamma$. Let $Q(x, y)= |y|^{\gamma}$. Let $u$ be an $\epsilon_0$-local minimizer of $J_Q(\cdot, B_2^{n}(0,0))$. Then, there is an integer $N_0=N_0(n, \gamma)$ such that for all
\begin{align*}
0< r \le c(n, \sup_{\partial B_2^n(0,0)}u_0, ||\nabla u_0||_{L^2(B_2^n(0, 0))}, \epsilon_0)
\end{align*}
such that $B_{r}^{n}(x, 0) \subset B_1^n(0,0)$, the ball $B_{r}^{n}(x, 0)$ does not have Property P with constant $N \ge N_0(n, \gamma)$.
\end{theorem}

\begin{remark}
Note that Property P and therefore Theorem \ref{t:main theorem 2} are geometric in nature, and not limited to density zero points.  Among other things, Theorem \ref{t:main theorem 2} implies (by Lemma \ref{l:technical heart}) that if the free boundary $\partial \{u>0 \}$ intersects $\Gamma$ at $(x, 0) \in B_2^n(0,0)$, then there is a constant $0<\alpha(n, \gamma)= (4N_0(n,\gamma))^{-1}$ such that for all sufficiently small scales no component of $\{u>0\} \cap B_r^n(x, 0)$  is contained in the cone set
\begin{align*}
    C_{\alpha, r}(x,0) := \{(x', y') \in B_r^n(x, 0): |y'| \le \alpha(n, \gamma) |x'| \}.
\end{align*}
Note that this is much stronger than simply eliminating density zero points and begins to describe the geometry of the positivity set $\{u>0 \}$ near $\Gamma$ beyond its infinitesimal properties.
\end{remark}

\subsection{Outline}
The remainder if the paper is broken up into 4 parts.  In Section 2, we state three essential lemmata: Lemma \ref{l: 1}, Lemma \ref{l: 2}, and Lemma \ref{l: 3}.  Assuming these lemmata, we are able to give the proofs of Theorem \ref{t:main theorem} and Theorem \ref{t:main theorem 2}.  The rest of the paper is dedicated to proving these three lemmata.  In Section 3, we present preliminary definitions and prior results .  Section 4 is the technical heart of the paper.  It is dedicated to proving Lemma \ref{l: 1} and is principally concerned with the local geometry of the positivity set $\{u>0\}$ near a cusp point. Section 5 is dedicated to proving Lemma \ref{l: 2} and Lemma \ref{l: 3}, which concern various growth conditions on $u$ in the neighborhoods constructed in Section 4.

\subsection{Acknowledgements}
The first author acknowledges the Center for Nonlinear Analysis at Carnegie Mellon University for its support.  Furthermore, the first author thanks Giovanni Leoni and Irene Fonseca for their invaluable generosity, patience, and guidance.  

\section{Proof of Theorem \ref{t:main theorem 2} and Theorem \ref{t:main theorem} }

The proof of Theorem \ref{t:main theorem 2} and Theorem \ref{t:main theorem} follows from three lemmata.  The first remark allows us to normalize the window in which we work.

\begin{definition}\label{d:rescalings}(\cite{AltCaffarelli81} Remark 3.1)
Let $0< \gamma$, and $u$ be an $\epsilon_0$-local minimizer of (\ref{e: AC functional}) in the class $\mathcal{K}_{u_0, B_2^n(0, 0)}$.  Then, for any $B^n_r(x, 0) \in B_2^n(0,0)$ we define the rescalings
\begin{align}\label{e:rescaling}
u_{(x,0), r}(x',y') := \frac{u(r(x',y') + (x, 0))}{r^{\gamma +1}}.
\end{align}
The function $u_{(x, 0), r}$ is an $\epsilon'$-local minimizer of (\ref{e: AC functional}) in the class $K_{u_{(x, 0), r}, B_{\frac{1}{r}}^n((0,0))},$ where
\begin{align*}
    \epsilon'= \epsilon_0 \max\{r^{-n}, r^{-(n-2\gamma)}\}.
\end{align*}
\end{definition}

Thus, for any $\epsilon_0$-minimizer of $J_{Q}(\cdot, \Omega)$ if $(x, 0) \in \Omega \cap \partial \{u>0\}$ we may reduce to considering $u_{(x, 0), r}$ a $1$-local minimizer of $J_{Q}(\cdot, B_2^n(0,0))$.  Thus, for the remainder of the paper, we argue for this case, assuming that $\Omega = B^n_2(0,0)$.

\begin{lemma}\label{l: 1}
If $u$ is a $\epsilon_0$-local minimizer of (\ref{e: AC functional}) for $\Omega = B^n_2(0,0)$ and $(x, 0) \subset B^n_2(0,0)$ satisfies $(x,0) \in \Sigma$, then for every $N \in \mathbb{N}$ there exists a radius
\begin{align*}
0< r(N) \le c(n, \sup_{\partial B_2^n(0,0)}u_0, ||\nabla u||_{L^2(B_2^n(0,0))}, \epsilon_0)
\end{align*}
and a ball $B^n_r(x, 0) \subset B^n_2(0,0)$ such that $B^n_r(x, 0)$ satisfies Property P with constant $N$.
\end{lemma}

\begin{definition}
We make an auxiliary definition of certain cylindrical neighborhoods. For any $(x', 0) \in \mathbb{R}^{k} \times \mathbb{R}^{n-k}$, we define
\begin{align*}
    A_j(x') := \left(B^k_{2j}(x') \setminus B^k_{2(j-1)}(x')\right) \times B^{n-k}_1(0).
\end{align*}
\end{definition}

\begin{remark}\label{r: A_j volume estimates}
For all $j \in \mathbb{N}$ and any $(x', 0) \in \mathbb{R}^{k} \times \mathbb{R}^{n-k}$,
\begin{align*}
  \mathcal{H}^n(A_j(x')) \approx_n j^{k-1}.
\end{align*}
\end{remark}

The next two lemmata are crucial growth estimates.  The first is a statement on the growth of harmonic functions.

\begin{lemma}\label{l: 2}
Let $u$ is an $\epsilon_0$-local minimizer of (\ref{e: AC functional}) for $\Omega = B^n_2(0,0)$.  Suppose that
\begin{align*}
0< r \le c(n, \sup_{\partial B_2^n(0,0)}u_0, ||\nabla u_0||_{L^2(B_2^n(0,0))}, \epsilon_0),
\end{align*}
$N \in \mathbb{N}$, and $B_{r}^n(x, 0) \subset B_2^n(0,0)$ satisfies Property P with constant $N$.  Then, if we let $\mathcal{O}' \subset \{u >0\}$ be the component guaranteed by Property P with constant N and $\mathcal{O} = \frac{N}{r}[\mathcal{O}'-(x, 0)]$, then for all integers $1 \le j \le N/2$,
\begin{align*}
    \int_{A_j(0) \cap \mathcal{O}} u^2_{(x, 0), \frac{r}{N}} dV \ge c(n, \gamma)j^{n-1.}
\end{align*}
\end{lemma}

The next statement is simple result of standard estimates on local minimizers of (\ref{e: AC functional}) and the growth rate of the cylindrical annuli $A_j$.

\begin{lemma}\label{l: 3}
Let $u$ is an $\epsilon_0$-local minimizer of (\ref{e: AC functional}) for $\Omega = B^n_2(0,0)$.  Suppose that
\begin{align*}
0< r \le c(n, \sup_{\partial B_2^n(0,0)}u_0, ||\nabla u_0||_{L^2(B_2^n(0,0))}, \epsilon_0),
\end{align*}
$N \in \mathbb{N}$, and $B_{r}^n(x, 0) \subset B_2^n(0,0)$ satisfies the Property P with constant $N$. Then, if we let $\mathcal{O}' \subset \{u >0\}$ be the component guaranteed by Property P with constant N and $\mathcal{O} = \frac{N}{r}[\mathcal{O}'-(x, 0)]$, then for all integers $1 \le j \le N/2$
\begin{align*}
    \int_{A_j(0) \cap \mathcal{O}} u_{(x,0),\frac{r}{N}}^2 dV \le C(n)j^{k-1}.
\end{align*}
\end{lemma}

With these lemmata, we prove Theorems \ref{t:main theorem}, \ref{t:main theorem 2} as follows.

\subsection{Proof of Theorem \ref{t:main theorem} and Theorem \ref{t:main theorem 2}} 

We note that by Lemma \ref{l: 1}, it suffices to prove Theorem \ref{t:main theorem 2}. Let $u$ be an $\epsilon_0$-local minimizer in the class $\mathcal{K}_{u_0, B_2^n(0,0)}$.  Suppose that 
\begin{align*}
0< r \le c(n, \sup_{\partial B_2^n(0,0)}u_0, ||\nabla u_0||_{L^2(B_2^n(0,0))}, \epsilon_0),
\end{align*}
$N \in \mathbb{N}$, and $B_{r}^n(x, 0) \subset B_2^n(0,0)$ satisfies Property P with constant $N$.  Let $u_{(x, 0), \frac{r}{N}}$ be as in (\ref{e:rescaling}) and $\mathcal{O}$ be the component of $\{u_{(x, 0), \frac{r}{N}}>0\}$ which is guaranteed by Lemma \ref{l: 2} and Lemma \ref{l: 3}.

We consider the ratio,
\begin{align}\label{e: ratio}
    \frac{\int_{A_j(0) \cap \mathcal{O}}u_{(x, 0), \frac{r}{N}}^2 dx}{\int_{A_j(0)} 1 dx}.
\end{align}
By Lemma \ref{l: 3} and Remark \ref{r: A_j volume estimates} we have that for all $1 \le j \le N/2,$
\begin{align*}
\frac{\int_{A_j(0) \cap \mathcal{O}}u_{(x, 0), \frac{r}{N}}^2 dx}{\int_{A_j(0)} 1 dx}\le \frac{C(n, \gamma) j^{k-1}}{C(n)j^{k-1}} \lesssim_{n, \gamma} 1.
\end{align*}

On the other hand, by Lemma \ref{l: 2} and Remark \ref{r: A_j volume estimates}, we have that 
\begin{align*}
\frac{\int_{A_j(0) \cap \mathcal{O}}u_{(x, 0), \frac{r}{N}}^2 dx}{\int_{A_j(0)} 1 dx}\gtrsim_{n, \gamma} j^{n-k}.
\end{align*}
for all integers $1 \le j \le N/2$. Therefore, if $N$ is sufficiently large (depending only upon $n, \gamma$) we may take $j$ sufficiently large and obtain a contradiction. \qed

\section{Preliminaries}

\begin{definition}
The techniques used in \cite{AltCaffarelli81} to establish the non-degeneracy of a local minimizer $u$ rely upon comparing $u$ with two other functions:
\begin{enumerate}
    \item The harmonic extension of $u$ in a ball $B_r^n(x, y).$
    \item The function $w = \min\{u, v\}$ in $B_r^n(0,0)$ for
\begin{align*}
v(x, y) = \left(\sup_{(x', y') \in B_{r \sqrt{s}}^n(0, 0)} \{u(x',y')\}\right) \max\left\{1 - \frac{|(x, y)|^{2-n} - r^{2-n}}{(sr)^{2-n} - r^{2-n}} , 0\right\}.
\end{align*} 
\end{enumerate}
Since $\norm{\nabla v}^2_{L^2(B^n_r(0,0))} \le C(s, n) \sup_{y \in B^n_{r \sqrt{s}}(0,0)} \{u^2(y)\} r^{n-2}$ and harmonic functions are energy minimizers, for every $\epsilon_0$-local minimizer $u$, there is a uniform scale 
\begin{align} \label{e:standard scale}
    r_0:= r_0(n, \sup_{\partial B_2^n(0,0)} u_0, \norm{\nabla u}_{L^2(B^n_2(0,0))}, \epsilon_0)
\end{align}
at which we can apply these arguments. We shall refer to this scale $r_0$ as the \textit{standard scale}.
\end{definition}

\begin{remark}
Let $0< \epsilon_0$, and let $n, k$ be integers such that $n \ge 2$ and $0 \le k \le n-1$. Let $\Gamma = \{(x, 0) \in \mathbb{R}^k\times \mathbb{R}^{n-k}\}$, and $0< \gamma$. Let $Q(x, y)= |y|^{\gamma}$. Let $u$ be an $\epsilon_0$-local minimizer of $J_Q(\cdot, B_2^{n}((0,0)))$. There is a \begin{align*}
0< c(n, \sup_{\partial B_2^n(0,0)}u_0, ||\nabla u||_{L^2(B_2^n(0,0))}, \epsilon_0)    
\end{align*}
such that for any $0< r \le c$ and any $B_r^n(x, 0) \subset B_2^n(0,0)$ the function $u_{(x,0),r}$ is a $1$-local minimizer in $\mathcal{K}_{u_{(0,0),r}, B_{\frac{2}{r}}^n(-\frac{x}{r},0)}$ and the standard scale $r_0$ for $u_{(0,0),r}$ is $r_0 = 1.$
\end{remark}

\begin{definition} For any continuous $u: \mathbb{R}^n \rightarrow \mathbb{R}$ and any $(x, y) \in \mathbb{R}^k \times \mathbb{R}^{n-k}$ we define the quantity
\begin{align}
    H(r, (x, y), u):= \int_{\partial B_r^n(x, y)}u^2 d\sigma.
\end{align}
\end{definition}

\begin{lemma}\label{H growth} For any function $u$ which is an $\epsilon_0$-local minimizer of (\ref{e: AC functional}) the following hold.  The function $u$ is continuous, and for any $B_r^n(x, y) \subset B_R^n(x, y) \subset \subset B^n_2(0,0)$
\begin{align}\label{H growth}
    H(R, (x, y), u) \ge R^{n-1} \left(\frac{1}{r^{n-1}}H(r, (x,y), u)\right).
\end{align}
\end{lemma}

The continuity of $u$ is proven in \cite{Mccurdy20} Lemma 3.14.  The inequality (\ref{H growth}) follows from the argument provided in \cite{HanLin_nodalsets}  (see Corollary 2.2.5 and Corollary 2.2.6) for harmonic functions. 


\begin{lemma}\label{c:local Lipschitz bound}(Local Lipschitz, \cite{Mccurdy20} Corollary 3.14)
Let $u$ be a $1$-local minimizer of $J_{Q}(\cdot, B_2^n(0,0))$ with standard scale $r_0 = 1.$  Assume that $(x,y) \in \partial \{u>0\}$.  Then, for all  $(x',y') \in \{u>0\} \cap B_1^n(x, y)$ 
\begin{align*}
|\nabla u(x', y')| \le C_1(n)\max\{\dist((x', y'), \partial \{u>0\}), |y'|\}^{\gamma}.
\end{align*}
\end{lemma}

\begin{lemma}\label{l:interior balls}(Interior Balls, \cite{Mccurdy20} Lemma 3.17)
Let $u$ be a $1$-local minimizer of $J_{Q}(\cdot, B_2^n(0,0))$ with standard scale $r_0 = 1.$  Let $0<r < 1$, and let $(x,y) \in \partial \{u>0\}$ satisfy $B_{r}^n(x,y) \subset B_2^n(0,0) \setminus \Gamma$. Then there exists a point $(x', y') \in \{u>0 \} \cap \partial B_{\frac{1}{2}r}^n(x,y)$ and a constant $0< c(n, Q_{\min, B_r^n(x,y)}) <\frac{1}{2}$ such that
\begin{align*}
    u \ge C_1(n) r Q_{\min, B_{\frac{1}{2}r}^n(x,y)}.
\end{align*}
in $B_{\frac{c}{2}r}^n(x', y')$.
\end{lemma}

\section{Proof of Lemma \ref{l: 1}}

In this section, we investigate the local geometry of the positivity set $\{u>0\}$ near a point in $\Sigma$. 

\begin{lemma}(Attenuation Radius)\label{r:attenuation radius}
Suppose that $u$ is an $\epsilon_0$-local minimizer of (\ref{e: AC functional}) in $B_2^n(0,0)$ for $Q(x,y) = |y|^\gamma$, as above. Let $(0,0) \in \Sigma.$ For any $0<\eta$ we can find a radius $0< r(\eta)< c(n, \sup_{\partial B_2^n(0,0)}u_0, ||\nabla u||_{L^2(B_2^n(0,0))})$ such that 
\begin{align}
    \{u_{(0,0), r}>0\} \cap (B_1^k(0)\times B^{n-k}_1(0)) \subset \{(x, y) \in \mathbb{R}^k \times \mathbb{R}^{n-k}: |y| \le \eta |x|\}. 
\end{align}
\end{lemma}

\begin{proof}
Suppose that the claim is false.  Then, there exists an $0< \eta$ for which there is a sequence of radii  $r_j \rightarrow 0$ such that there exists a $(x_j, y_j) \in \partial \{u_{(0,0), r}>0\} \cap B_1^n(0, 0)$ for which $|y_j| \ge \eta|x_j|$. For sufficiently small $0<r < c(n, \sup_{\partial B_2^n(0,0)}u_0, ||\nabla u||_{L^2(B_2^n(0,0))}, \epsilon_0)    
$, then we may apply Corollary \ref{c:local Lipschitz bound} and Lemma \ref{l:interior balls},
\begin{align*}
    \frac{\mathcal{H}^n(B_1^n(0,0) \cap \{u_{(0,0), r_j}>0\})}{\omega_n}> c>0.
\end{align*}
This contradicts $(0, 0) \in \Sigma.$
\end{proof}

\begin{lemma}\label{l:technical heart}(\cite{McCurdy21} Lemma 3.4)
Let $n \ge 2$, and let $0 \le k \le n-1$.  Let $(x, y) \in \mathbb{R}^{k}\times \mathbb{R}^{n-k}.$  Let $0< \gamma$ and $Q(x, y) = |y|^\gamma$.  Let $u$ be an $\epsilon_0$-local minimizer of (\ref{e: AC functional}). Suppose that $(0, 0) \in \Sigma$ and let $\mathcal{O}$ be any component of $\{u>0\}$ such that $(0, 0) \in \overline{\mathcal{O}}$.

For any $1<N< \infty$, we may find a radius $0< \rho(N) \le r(1/4N)$ such that the rescaling $u_{(0, 0), \rho}$ satisfies the following conditions
\begin{itemize}
    \item [i.] There exists a large radius $1<<N_0<\infty$ satisfying $4N<N_0$ such that
\begin{align}\label{e: standard window upper bound}
    \max \{|y| : (x, y) \in \left(\partial B^k_{N_0}(0) \times [-2,2]^{n-k}\right) \cap \partial \mathcal{O}_{(0, 0), \rho}\} &= 1.
\end{align}
    \item [ii.] For the radius $N_0 - N$
\begin{align}
    \max \{|y| : (x, y) \in \left(\partial B^k_{N_0-N}(0) \times [-2,2]^{n-k}\right) \cap \partial \mathcal{O}_{(0, 0), \rho}\} &\ge \frac{1}{2}.
\end{align}
\end{itemize}
\end{lemma}

In the next lemma, we derive upper and lower bounds on both the supremum of the distance away from $\Gamma$ and the value of $u_{(0,0),\rho}$ inside a specified region of $\mathcal{0}$. To ease notation we make the following definitions.
Set
\begin{align*}
    \text{Height}(r, \mathcal{O}_{(0,0), \rho}, R):= & \sup \{|y| : (x, y) \in \left(\partial B^k_{\rho}(0) \times B^{n-k}_{R}(0)\right) \cap \mathcal{O}_{(0,0), \rho}\}\\
    M(r, \mathcal{O}_{(0,0), \rho}, R) := & \sup \{u_{(0,0), \rho}(x, y) : (x, y) \in \left(\partial B^k_{\rho}(0) \times B^{n-k}_{R}(0) \right) \cap \mathcal{O}_{(0,0), \rho}\}.   
\end{align*}
\begin{lemma}\label{l:height bound}(Height bound)
Let $u$ is a local minimizer of (\ref{e: AC functional}). Suppose that $(0,0) \in \Sigma$. Let $\mathcal{O}$ be any component of $\{u>0 \}$ such that $(0,0) \in \overline{\mathcal{O}}.$ Let $2^\gamma + 1/4 < N \in \mathbb{N}$ be fixed.  Let $0< \rho =\rho(N)$ as in Lemma \ref{l:technical heart}. 
Then,
\begin{align}\nonumber
    1/4 & < \inf_{r \in [N_0 - N + \frac{1}{4}, N_0]} \text{\emph{Height}}(r, \mathcal{O}_{(0,0), \rho}, 4N) \\ &  \le \sup_{r \in [N_0 - N, N_0-2^{\gamma}]} \text{\emph{Height}}(r, \mathcal{O}_{(0,0), \rho}, 4N)< 2^{1+\gamma}< 4N
\end{align}
and
\begin{align}\nonumber
    C_1(n)\left(\frac{1}{4}\right)^{1+\gamma} & \le \min_{r \in [N_0 - N + \frac{1}{4}, N_0]}M(r, \mathcal{O}_{(0,0), \rho}, 4N)\\ & \le \max_{r \in [N_0 - N + \frac{1}{4}, N_0]}M(r, \mathcal{O}_{(0,0), \rho}, 4N) \le C_1(n).
\end{align}
\end{lemma}

\begin{proof}
Let $N \in \mathbb{N}$ be given and consider $u_{(0,0), \rho}$ and $\mathcal{O}_{(0,0), \rho}$ as in Lemma \ref{l:technical heart}. Because $u_{(0,0), \rho}$ is harmonic in $\mathcal{O}_{(0, 0), \rho}$ and vanishes on $\partial \mathcal{O}_{(0,0), \rho}$, the Maximum Principle implies $M(r_1, \mathcal{O}_{(0,0), \rho}, 4N) \le M(r_2, \mathcal{O}_{(0,0), \rho}, 4N)$ for all $0 \le r_1 < r_2 \le N_0$. 

Now, we observe that Lemma \ref{l:interior balls} implies that if $\text{Height}(r, \mathcal{O}_{(0,0), \rho}, 4N) = c_0$ then 
\begin{align*}
    C_1(n)\left(\frac{c_0}{2}\right)^{1+\gamma} \le  M\left(r + \frac{c_0}{2}\right).
\end{align*}
Moreover, by Lemma \ref{c:local Lipschitz bound}
\begin{align*}
    M(r, \mathcal{O}_{(0,0), \rho}, 4N) \le C_1(n)c_0^{1+\gamma}.
\end{align*}
Therefore, since $\text{Height}(N_0, \mathcal{O}_{(0,0), \rho}, 4N) = 1$, we have 
\begin{align*}
    M(N_0, \mathcal{O}_{(0,0), \rho}, 4N) \le C_1(n)
\end{align*}
which implies that $\text{Height}(r, \mathcal{O}_{(0,0), \rho}, 4N) \le 2^{1+\gamma}$ for all $0<r \le N_0-2^\gamma.$

Similarly, since $\text{Height}(N_0 - N, \mathcal{O}_{(0,0), \rho}, 4N) \ge \frac{1}{2}$, we have 
\begin{align*}
    M(N_0-N+ \frac{1}{4}, \mathcal{O}_{(0,0), \rho}, 4N) \ge C_1(n)\left(\frac{1}{4}\right)^{1+\gamma}
\end{align*}
which implies that $\text{Height}(r, \mathcal{O}_{(0,0), \rho}, R) \ge 1/4$ for all $N_0-N+1/4 \le r \le N_0$. This proves the lemma.
\end{proof}

\begin{corollary}
\label{c: height bound}
Let $n \ge 2$, and let $0 \le k \le n-1$. Let $0< \gamma$ and $Q(x, y) = |y|^\gamma$.  Let $u$ be an $\epsilon_0$-local minimizer of (\ref{e: AC functional}). Suppose that $B_r^n(x, 0) \subset B_2^n(0,0)$ satisfies Property P with constant $N$.  Let $\mathcal{O}'$ be the component of $\{u>0\}$ guaranteed by Property P, and let $\mathcal{O} = \frac{N}{r}[\mathcal{O}-(x, 0)]$. By the rescaling in Definition \ref{d:rescalings}, then
\begin{align}\nonumber
    \frac{1}{4\cdot 2^{1+\gamma}} & < \inf_{r \in [0, N]} \text{\emph{Height}}(r, \mathcal{O}, 2) \le \sup_{r \in [0, N-2^{\gamma}]} \text{\emph{Height}}(r, \mathcal{O}, 2)< 1
\end{align}

Moreover, if $0<r \le c(n, \sup_{\partial B_2^n(0,0)}u_0, ||\nabla u||_{L^2(B_2^n(0,0))}, \epsilon_0)$ then the proof of Lemma \ref{l:height bound} shows that
\begin{align}\label{Prop P u value bounds}
    C_1(n)\left(\frac{1}{4}\right)^{1+\gamma} \le \min_{r \in [0, N]}M(r, \mathcal{O}, 2) \le \max_{r \in [0, N]}M(r, \mathcal{O}, 2) \le C_1(n).
\end{align}
\end{corollary}

\subsection{Proof of Lemma \ref{l: 1}}. To prove Lemma \ref{l: 1}, we let $N \in \mathbb{N}$ be given.  We then consider $\tilde{N} \in \mathbb{N}$ such that $2^{2+\gamma}N + 2^{\gamma} + 1/4 \le \tilde{N}$.  Choosing this $\tilde{N}$ for Lemma \ref{l:technical heart} and Lemma \ref{l:height bound}, we obtain an $0< \rho(\tilde{N})$ and $u_{(0,0), \rho}$ and $\mathcal{O}_{(0,0), \rho}$ such that for all $N_0 - \tilde{N} + 1/4 \le r \le N_0 -2^{\gamma}$
\begin{align*}
    1/4 \le \text{Height}(r, \mathcal{O}_{(0,0), \rho}, 4N)\le 2^{1+\gamma}.
\end{align*}
Let $(x', y) \in \mathbb{R}^k \times \mathbb{R}^{n-k} \cap \overline{\mathcal{O}_{(0,0), \rho}}$ such that $|x'| = \frac{1}{2}((N_0 - 2^{\gamma}) + (N_0- N + 1/4))$ and $y$ is chosen such that $|y| = \text{Height}(|x'|, \mathcal{O}_{(0,0), \rho}, 4N)$.

The rescaling $u_{(0,0), \frac{\rho}{2^{1+\gamma}}}$ with $(x'/2^{1+\gamma}, 0)$ satisfies the conditions of Lemma \ref{l: 1}.

\section{Proof of Lemma \ref{l: 2} and Lemma \ref{l: 3}}

In the section, we present the proofs of Lemma \ref{l: 2} and Lemma \ref{l: 3}.  First, we prove an auxiliary lemma providing a lower bound on the thickness of spherical shells that which can be contained with the annuli $A_j(x)$.

\begin{lemma}
Let $(x,0) \in \mathbb{R}^{k} \times \mathbb{R}^{n-k}$.  Then, there exists a constant $0< c \le 1$ such that for all $j \in \mathbb{N}$
\begin{align*}
    0< c \le \mathcal{H}^1\left(\{0<r: \partial B_r^n(x, y) \cap B_1^n(\Gamma) \subset A_j(x)\}\right).
\end{align*}
\end{lemma}

\begin{proof}
A point $(x',y')$ is contained in $A_j(x)$ if $2(j-1)\le\sqrt{(x_1')^2+...+(x_k')^2}\le 2j$ and $0\le\sqrt{(y_{1}')^2+...+(y_{n-k}')^2}\le 1$.  Thus $\partial B_r^n(x, y) \cap B_1^n(\Gamma)\subset A_j(x)$ when $(2j-2)^2+1\le r^2\le (2j)^2$.
This means we can take $0<c=\min_{j\in\mathbb{N}}|2j - \sqrt{(2(j-1))^2+1}|=1$
\end{proof}

\subsection{Proof of Lemma \ref{l: 2}}

Let $B_r^n(x, 0)$, $\mathcal{O}$, and $u_{(x, 0), r}$ be as in Lemma \ref{l: 2}. For ease of notation we define $R_j=2j$ and $r_j=\sqrt{(2(j-1))^2+1}$
\begin{align*}
\int_{A_j(0)\cap\mathcal{O}}u_{(x, 0), r}^2dV
 &\ge\int_{r_j}^{R_j}\int_{\partial B_\rho^n(0,0)\cap\mathcal{O}}u_{(x, 0), r}^2 d\sigma d\rho\\
 &=\int_{r_j}^{R_j} H(\rho,(0, 0), u_{(x, 0), r})d\rho\\
 &\ge \int_{r_j}^{R_j}\frac{\rho^{n-1}}{r_1^{n-1}}H(r_1,(0,0), u_{(x, 0), r})d\rho\\
 &\ge C(n)\frac{j^{n-1}}{r_1^{n-1}}H(r_1, (0, 0), u_{(x, 0), r}).
\end{align*}
By Corollary \ref{Prop P u value bounds}, $M(0)= \sup \{u_{(x, 0), r}(0, y): |y|\le 1\} \ge \frac{C_1(n)}{4^{1+ \gamma}}$.  Hence, by the Maximum principle and the continuity of $u_{(x, 0), r}$ there must be a point $(v,w) \in \partial B_{r_1}^n(0,0)$ such that
\begin{align*}
    u_{(x, 0), r}(v,w) = \max_{\partial B_{r_1}^n(0,0)} u_{(x,0), r} \ge \frac{C_1(n)}{4^{1+ \gamma}}.
\end{align*}
By Lemma \ref{c:local Lipschitz bound}, $B_{\frac{1}{4^{1+\gamma}}}^n(v,w) \subset \{u_{(x,0), r}>0 \}$ and  
\begin{align*}
\min_{B_{\frac{1}{2\cdot 4^{1+\gamma}}}^n(v,w)} u_{(x,0), r} \ge \frac{C_1(n)}{2\cdot 4^{1+ \gamma}}\\
     \frac{1}{r_1^{n-1}}H(r_1, (0, 0), u_{(x, 0), r}) \ge c(n, \gamma).
\end{align*}
This concludes the proof of Lemma \ref{l: 2}.

\subsection{Proof of Lemma \ref{l: 3}}
By Corollary \ref{c: height bound} we get that $u_{(x,0),r}\le C_1(n)$.  Then by 
\[\int_{A_j(0)\cap\mathcal{O}}u_{(x,0), r}^2dV\lesssim_{n}\int_{A_j(0)\cap\mathcal{O}}1dV \le \int_{A_j(0)} 1 dV\le C(n)j^{k-1}.\] 
This concludes the proof of Lemma \ref{l: 3}.

\bibliography{references}
\bibliographystyle{amsalpha}

\end{document}